\documentclass[12pt,reqno]{amsart}
\usepackage{amsfonts,amssymb,amsmath}

\DeclareMathOperator{\bmax}{\mathbf{max}}
\DeclareMathOperator{\bmin}{\mathbf{min}}
\DeclareMathOperator{\Filters}{\mathcal{F}}
\DeclareMathOperator{\Ant}{\mathfrak{A}}
\DeclareMathOperator{\Antb}{\mathfrak{B}}
\DeclareMathOperator{\Antc}{\mathfrak{C}}

\DeclareMathOperator{\bmap}{\mathfrak{b}}
\DeclareMathOperator{\cmap}{\mathfrak{c}}

\newcommand{\Pa}{{P^{\mathrm a}}}
\newcommand{\Qa}{{Q^{\mathrm a}}}

\newtheorem{thm}{Theorem}
\newtheorem{defn}[thm]{Definition}

\newtheorem{lem}[thm]{Lemma}
\newtheorem{prop}[thm]{Proposition}

\numberwithin{equation}{section} \numberwithin{thm}{section}

\title{On blockers in bounded posets}

\date{}
\author{Andrey O. Matveev}
\address{Data-Center Company, RU-620034, P.O.~Box~5, Ekaterinburg,
Russia} \email{aomatveev@dc.ru}

\begin{document}

\begin{abstract}
Antichains of a finite bounded poset are assigned antichains
playing a role analogous to that played by blockers in the Boolean
lattice of all subsets of a finite set. Some properties of
lattices of generalized blockers are discussed.
\end{abstract}

\maketitle

\section{Introduction}

Blocking sets for finite families of finite sets are important
objects of discrete mathematics
(see~\cite[Chapter~8]{GLS},~\cite{F}).

A set $H$ is called a {\em blocking set\/} for a nonempty family
$\mathcal{G}=\{G_1,\ldots,G_m\}$ of nonempty subsets of a finite set if for
each $k\in\{1,\ldots,m\}$ we have $|H\cap G_k|\geq 1$. The {\em blocker\/} of
$\mathcal{G}$ is the family of all inclusion-wise minimal blocking sets for
$\mathcal{G}$.

A family of subsets of a finite set is called a {\em clutter\/}
(or a {\em Sperner family}) if no set from it contains another. If
the family is empty or if it consists of only one subset,
$\{\emptyset\}$, then the corresponding clutter is called {\em
trivial}.

The concepts of {\em blocker map\/} and {\em complementary map\/} on
clutters~\cite{CFM} made it possible to clarify the relationship between
specific families of sets, arising from the matroid theory, and maps on them.
The blocker map, that assigns the blocker to a clutter, is defined on all
clutters, including trivial clutters.

The following property~\cite{EF,L} is basic: for a clutter
$\mathcal{G}$, the blocker of its blocker coincides with
$\mathcal{G}$.

We show that the concepts of blocking set and blocker can be
extended when passing from discussing clutters, considered as
antichains of the Boolean lattice of all subsets of a finite set,
to exploring antichains of arbitrary finite bounded posets (a
poset $P$ is called {\em bounded\/} if it has a unique minimal
element, denoted $\hat{0}_P$, and a unique maximal element,
denoted $\hat{1}_P$).

In Section~\ref{Intersectersandcomplementers}, the notion of
intersecter plays a role analogous to that played by the notion of
blocking set in the Boolean lattice of all subsets of a finite
set. In Section~\ref{IntersectersinCartesianproductsofposets}, we
explore the structure of subposets of intersecters in Cartesian
products of posets. In
Section~\ref{Blockermapandcomplementarymap}, some properties of
the blocker map and complementary map are shortly discussed. In
Section~\ref{Latticeofblockers}, the structure of lattices of
generalized blockers is reviewed.

\section{Intersecters and complementers}
\label{Intersectersandcomplementers}

We refer the reader to~\cite[Chapter~3]{St} for basic information
and terminology in the theory of posets.

For a poset $Q$, $\Qa$ denotes its atom set; $\bmin Q$ and $\bmax
Q$ denote the sets of all minimal elements and all maximal
elements of $Q$, respectively; $\mathfrak{I}_Q(X)$ and
$\mathfrak{F}_Q(X)$ denote the order ideal and order filter of $Q$
generated by a subset $X\subseteq Q$, respectively. If $x,y$ are
elements of $Q$ and $x<y$ (or $x\leq y$) then we write $x<_Q y$
(or $x\leq_Q y$). In a similar way, we denote by $\vee_Q$ the
operation of join in a join-semilattice $Q$, and we denote by
$\wedge_Q$ the operation of meet in a meet-semilattice $Q$. We use
$\times$ to denote the operation of Cartesian product of posets.

For a finite family $\mathcal{G}$ of finite sets, its conventional
blocker is denoted by $\mathcal{B}(\mathcal{G})$.

Throughout $P$ stands for a finite bounded poset with $|P|>1$. We
start with extending of the concept of blocking set.

\begin{defn}
\label{Intersectersandcomplementers1} Let $A$ be a subset of $P$.

\begin{itemize}
\item
If $A\neq \emptyset$ and $A\neq\{\hat{0}_P\}$ then an element
$b\in P$ is an {\em intersecter for $A$ in $P$\/} if for every
$a\in A-\{\hat{0}_P\}$, we have
\begin{equation}
|\mathfrak{I}_P(b)\cap\mathfrak{I}_P(a)\cap \Pa|\geq 1\ .
\end{equation}
\item
If $A=\{\hat{0}_P\}$ then $A$ has no intersecters in $P$.
\item
If $A=\emptyset$ then every element of $P$ is an {\em intersecter for $A$ in
$P$}.

\item
Every non-intersecter for $A$ in $P$ is a {\em complementer for $A$ in $P$}.
\end{itemize}
\end{defn}

Let $\mathcal{L}$ denote a finite Boolean lattice. If $A$ is a nonempty subset
of the poset $\mathcal{L}-\{\hat{0}_{\mathcal{L}}\}$ then an element
$b\in\mathcal{L}$ is an intersecter for $A$ in $\mathcal{L}$ if and only if\/
$\mathfrak{I}_{\mathcal{L}}(b)\cap\mathcal{L}^{\mathrm a}$ is a blocking set
for the family $\{\mathfrak{I}_{\mathcal{L}}(a)\cap\mathcal{L}^{\mathrm a}:a\in
A\}$.

We denote by $\mathbf{I}(P,A)$ and $\mathbf{C}(P,A)$ the sets of
all intersecters and all complementers for $A$ in $P$,
respectively. We consider the sets $\mathbf{I}(P,A)$ and
$\mathbf{C}(P,A)$ as subposets of the poset $P$. For a one-element
set $\{a\}$ we write $\mathbf{I}(P,a)$ instead of
$\mathbf{I}(P,\{a\})$ and $\mathbf{C}(P,a)$ instead of
$\mathbf{C}(P,\{a\})$.

We have the partition $\mathbf{I}(P,A)\ \dot{\cup}\
\mathbf{C}(P,A)=P$. For a nonempty subset $A\subseteq
P-\{\hat{0}_P\}$, the subposets of all its intersecters and
complementers are nonempty; indeed, we have
$\mathbf{I}(P,A)\ni\hat{1}_P$ and $\mathbf{C}(P,A)\ni\hat{0}_P$.
It follows from Definition~\ref{Intersectersandcomplementers1}
that for such a subset $A$, we have
\begin{equation}
\mathbf{I}(P,A)=\mathbf{I}(P,\bmin A)\ ,\
\mathbf{C}(P,A)=\mathbf{C}(P,\bmin A)\ ,
\end{equation}
therefore, in most cases, we may restrict ourselves to considering intersecters
and complementers for antichains; further,
\begin{equation}
\label{Intersectersandcomplementers2} \mathbf{I}(P,A)=\bigcap_{a\in
A}\mathbf{I}(P,a)\ ,\ \mathbf{C}(P,A)=\bigcup_{a\in A}\mathbf{C}(P,a)\ .
\end{equation}

For all antichains (including the empty antichain) $A_1,A_2$ of $P$ with
$\mathfrak{F}_P(A_1)\subseteq\mathfrak{F}_P(A_2)$, we have
\begin{equation}
\label{Intersectersandcomplementers3}
\mathbf{I}(P,A_1)\supseteq\mathbf{I}(P,A_2)\ ,\
\mathbf{C}(P,A_1)\subseteq\mathbf{C}(P,A_2)\ .
\end{equation}

Clearly, the subposet $\mathbf{I}(P,a)$ of all intersecters for an
element $a\in P$ is the order filter
$\mathfrak{F}_P\bigl(\mathfrak{I}_P(a)\cap\Pa\bigr)$, hence, in
view of~(\ref{Intersectersandcomplementers2}),
equality~(\ref{eq:2006-1}) in the following lemma holds.

\begin{lem}
\label{Intersectersandcomplementers4} Let $A$ be a nonempty subset of $P-
\{\hat{0}_P\}$. The subposet of all intersecters for $A$ in $P$ is determined
by the following equivalent equalities:

\begin{equation}
\label{eq:2006-1} \mathbf{I}(P,A)=\bigcap_{a\in A}\mathfrak{F}_P
\bigl(\mathfrak{I}_P(a)\cap\Pa\bigr)\ ,
\end{equation}

\begin{equation}
\label{eq:2006-2}
\mathbf{I}(P,A)=
\bigcup_{E\in\mathcal{B}(\{\mathfrak{I}_P(a)\cap\Pa:a\in A\})}
\quad\bigcap_{e\in E}\mathfrak{F}_P(e)\ .
\end{equation}
\end{lem}

\begin{proof}
To prove~(\ref{eq:2006-2}), note that the inclusion
\begin{equation}
\mathbf{I}(P,A)\supseteq
\bigcup_{E\in\mathcal{B}(\{\mathfrak{I}_P(a)\cap\Pa:a\in A\})}\quad
\bigcap_{e\in E}\mathfrak{F}_P(e)
\end{equation}
follows from the definition of intersecters.

We are left with proving the inclusion
\begin{equation}
\mathbf{I}(P,A)\subseteq
\bigcup_{E\in\mathcal{B}(\{\mathfrak{I}_P(a)\cap\Pa:a\in A\})}\quad
\bigcap_{e\in E}\mathfrak{F}_P(e)\ .
\end{equation}
Assume that it does not hold, and consider such an intersecter $b$
for $A$ that
$b\not\in\bigcup_{E\in\mathcal{B}(\{\mathfrak{I}_P(a)\cap\Pa:a\in
A\})} \bigcap_{e\in E}\mathfrak{F}_P(e)$. In this case, the
inclusion $b\in\bigcap_{e\in E}\mathfrak{F}_P(e)$ holds not for
all sets $E$ from the family
$\mathcal{B}(\{\mathfrak{I}_P(a)\cap\Pa:a\in A\})$, hence there
exists such an element $a\in A$ that
$|\mathfrak{I}_P(b)\cap\mathfrak{I}_P(a)\cap\Pa|=0$. Therefore $b$
is not an intersecter for $A$, but this contradicts our choice of
$b$. Hence, (\ref{eq:2006-2}) holds.
\end{proof}

Thus, for every antichain $A$ of the poset $P$, the subposet of
all intersecters for $A$ in $P$ is an order filter of $P$, that
is,
$\mathbf{I}(P,A)=\mathfrak{F}_P\bigl(\bmin\mathbf{I}(P,A)\bigr)$.
As a consequence, the subposet $\mathbf{C}(P,A)$ of all
complementers for $A$ in $P$ is the order ideal
$\mathfrak{I}_P\bigl(\bmax\mathbf{C}(P,A)\bigr)$.

If $A$ is a subset of the poset $P$ then we call the antichain
$\bmin\mathbf{I}(P,A)$ the {\em blocker of $A$ in $P$}. We call
elements of the blocker $\bmin\mathbf{I}(P,A)$ {\em minimal
intersecters for $A$ in $P$}, and we call elements of the
antichain $\bmax\mathbf{C}(P,A)$ {\em maximal complementers for
$A$ in $P$}.

The images of intersecters under suitable order-preserving maps are also
intersecters:

\begin{prop}
Let $P_1$ and $P_2$ be disjoint finite bounded posets with $|P_1|,|P_2|>1$. Let
$\psi:P_1\rightarrow P_2$ be an order-preserving map such that
\begin{equation}
\label{Intersectersandcomplementers7}
\psi(\hat{0}_{P_1})=\hat{0}_{P_2}\ ,\ \ \ \ \
\psi(x_1)>_{P_2}\hat{0}_{P_2},\ \ \forall x_1>_{P_1}\hat{0}_{P_1}\
.
\end{equation}
For every subset $A_1$ of $P_1$
\begin{equation}
\psi\bigl(\mathbf{I}(P_1,A_1)\bigr)\subseteq
\mathbf{I}\bigl(P_2,\psi(A_1)\bigr)\ .
\end{equation}
\end{prop}

\begin{proof}
There is nothing to prove for $A_1=\emptyset\subset P$ or
$A_1=\{\hat{0}_{P_1}\}$. So suppose that $A_1\neq\emptyset\subset
P$ and $A_1\neq\{\hat{0}_{P_1}\}$. Let $b_1$ be an intersecter for
$A_1$. According to
Definition~\ref{Intersectersandcomplementers1}, for all $a_1\in
A_1$, $a_1>_{P_1}\hat{0}_{P_1}$, we have
$|\mathfrak{I}_{P_1}(b_1)\cap\mathfrak{I}_{P_1}(a_1)\cap
P_1{}^{\mathrm a}|\geq 1$, and in view
of~(\ref{Intersectersandcomplementers7}), for every atom $z_1\in
\mathfrak{I}_{P_1}(b_1)\cap\mathfrak{I}_{P_1}(a_1)\cap
P_1{}^{\mathrm a}$ we have the inclusion
\begin{equation}
\mathfrak{I}_{P_2}\bigl(\psi(z_1)\bigr)\cap P_2{}^{\mathrm
a}\subseteq\mathfrak{I}_{P_2} \bigl(\psi(a_1)\bigr)\cap P_2{}^{\mathrm a}\ ,
\end{equation}
the left-hand part of which is nonempty. Hence, for all $a_2\in\psi(A_1)$ the
inclusion $b_1\in\mathbf{I}(P,A_1)$ implies that
\begin{equation}
\left|\mathfrak{I}_{P_2}\bigl(\psi(b_1)\bigr)
\cap\mathfrak{I}_{P_2}\bigl(\psi(a_1)\bigr)\cap P_2{}^{\mathrm a}\right|\geq 1\
.
\end{equation}
This means that $\psi(b_1)\in\mathbf{I}\bigl(P_2,\psi(A_1)\bigr)$ and completes
the proof.
\end{proof}

\section{Intersecters in Cartesian products of posets}
\label{IntersectersinCartesianproductsofposets}

In this section, we study the structure of subposets of
intersecters in Cartesian products of two finite posets.

\begin{prop}
Let $P_1$ and $P_2$ be disjoint
finite bounded posets with $|P_1|,|P_2|>2$. Let $Q$ denote the poset
\begin{equation}
(P_1-\{\hat{0}_{P_1},\hat{1}_{P_1}\})\times
(P_2-\{\hat{0}_{P_2},\hat{1}_{P_2}\})\ \dot{\cup}\ \{\hat{0}_Q,\hat{1}_Q\}\ ,
\end{equation}
where $\hat{0}_Q$ and $\hat{1}_Q$ are the adjoint new least and greatest
elements. Let $A$ be a nonempty subset of the poset
$Q-\{\hat{0}_Q,\hat{1}_Q\}$, and let $A\!\!\downharpoonright_{P_1}$ and
$A\!\!\downharpoonright_{P_2}$ denote the subsets $\{a_1\in P_1:(a_1;a_2)\in
A\}$ and $\{a_2\in P_2:(a_1;a_2)\in A\}$, respectively.

\begin{enumerate}
\renewcommand{\theenumi}{\rm{\roman{enumi}}}

\item
If $\bmin\mathbf{I}(P_1,
A\!\!\downharpoonright_{P_1})=\{\hat{1}_{P_1}\}$ or
$\bmin\mathbf{I}(P_2,
A\!\!\downharpoonright_{P_2})=\{\hat{1}_{P_2}\}$, then
\begin{equation}
\mathbf{I}(Q,A)=\bmin\mathbf{I}(Q,A) =\{\hat{1}_Q\}\ .
\end{equation}

\item
If $\bmin\mathbf{I}(P_1,
A\!\!\downharpoonright_{P_1})\neq\{\hat{1}_{P_1}\}$ and
$\bmin\mathbf{I}(P_2,
A\!\!\downharpoonright_{P_2})\neq\{\hat{1}_{P_2}\}$, then
\begin{equation}
\mathbf{I}(Q,A)= \left(\mathbf{I}(P_1,
A\!\!\downharpoonright_{P_1})-\{\hat{1}_{P_1}\}\right) \times \left(\mathbf{I}(
P_2, A\!\!\downharpoonright_{P_2})-\{\hat{1}_{P_2}\}\right)\ \dot{\cup}\
\{\hat{1}_Q\}\ ,
\end{equation}
and $\bmin\mathbf{I}(Q,A)= \bmin\mathbf{I}(P_1,
A\!\!\downharpoonright_{P_1})\times \bmin\mathbf{I}(P_2,
A\!\!\downharpoonright_{P_2})$.
\end{enumerate}
\end{prop}

\begin{proof}
The atom set $\Qa$ of the poset $Q$ is $P_1{}^{\mathrm a}\times
P_2{}^{\mathrm a}$, therefore, by (\ref{eq:2006-1}), the subposet
of intersecters for $A$ in $Q$ is
\begin{align}
\mathbf{I}(Q,A)&= \bigcap_{a=(a_1;a_2)\in A}\mathfrak{F}_{Q} \biggl(
\Bigl(\mathfrak{I}_{P_1}(a_1)\times\mathfrak{I}_{P_2}(a_2)\Bigr)\cap \Bigl(
P_1{}^{\mathrm a}\times P_2{}^{\mathrm a}\Bigr) \biggr)
\\
&= \bigl( \mathbf{I}(P_1,A\!\!\downharpoonright_{P_1}) -\{\hat{1}_{P_1}\}
\bigr) \times \bigl( \mathbf{I}(P_2,A\!\!\downharpoonright_{P_2})
-\{\hat{1}_{P_2}\} \bigr)\ \dot{\cup}\ \{\hat{1}_Q\}\ ,
\end{align}
and the statement follows.
\end{proof}

\begin{prop}
Let $P_1$ and $P_2$ be disjoint finite bounded posets with
$|P_1|,|P_2|>1$. Let $Q$ denote the poset $P_1\times P_2$, and let
$A$ be a nonempty subset of the poset $Q-\{\hat{0}_Q\}$. Then

\begin{equation}
\mathbf{I}(Q,A)=\bigcap_{(a_1;a_2)\in
A}\Bigl(\bigl(P_1\times\mathbf{I}(P_2,a_2)\bigr)\ \cup\
\bigl(\mathbf{I}(P_1,a_1)\times P_2\bigr)\Bigr)\ .
\end{equation}
\end{prop}

\begin{proof}
Since the atom set $\Qa$ of the poset $Q$ is
$\bigl(\{\hat{0}_1\}\times P_2{}^{\mathrm a}
\bigr)\dot{\cup}\bigl( P_1{}^{\mathrm
a}\times\{\hat{0}_2\}\bigr)$, we have, according to
equality~(\ref{eq:2006-1}),
\begin{equation}
\begin{split}
\mathbf{I}(Q,A)&= \bigcap\limits_{ (a_1;a_2)\in A}\mathfrak{F}_{Q} \biggl(
\Bigl(\mathfrak{I}_{P_1}(a_1)\times\mathfrak{I}_{P_2}(a_2)\Bigr)\cap \Bigl(
\bigl(\{\hat{0}_1\}\times P_2{}^{\mathrm a}\bigr)\ \dot{\cup}\ \bigl(
P_1{}^{\mathrm a}\times\{\hat{0}_2\}\bigr) \Bigr) \biggr)\\ &=
\bigcap_{(a_1;a_2)\in
A}\mathfrak{F}_Q\Bigl(\bigl(\{\hat{0}_1\}\times(\mathfrak{I}_{P_2}(a_2)\cap
P_2{}^{\mathrm a})\bigr)\ \dot\cup\ \bigl((\mathfrak{I}_{P_1}(a_1)\cap
P_1{}^{\mathrm a})\times \{\hat{0}_2\}\bigr)\Bigr)\ ,
\end{split}
\end{equation}
and the statement follows.
\end{proof}

\section{Blocker map and complementary map}
\label{Blockermapandcomplementarymap}

Let $\Filters(P)$ denote the distributive lattice of all order
filters (partially ordered by inclusion) of $P$, and let $\Ant(P)$
denote the lattice of all antichains of $P$. For antichains
$A_1,A_2\in\Ant(P)$, we set
\begin{equation}
A_1\leq_{\Ant(P)}A_2\text{\ \ iff\ \ }
\mathfrak{F}_P(A_1)\subseteq\mathfrak{F}_P(A_2)\ ;
\end{equation}
in other words, we make use of the isomorphism
$\Filters(P)\rightarrow\Ant(P)$: $F\mapsto\bmin F$. We call the
least element $\hat{0}_{\Ant(P)}=\emptyset\subset P$ and greatest
element $\hat{1}_{\Ant(P)}=\{\hat{0}_P\}$ of the lattice
$\mathfrak{A}(P)$ the {\em trivial antichains\/} of $P$. They are
counterparts of trivial clutters.

Recall (cf.~\cite{GK}) that for $A_1,A_2\in\Ant(P)$,
\begin{equation}
A_1\vee_{\Ant(P)}A_2=\bmin(A_1\cup A_2)\ ,\ \ \
A_1\wedge_{\Ant(P)}A_2=\bmin\bigl(\mathfrak{F}_P(A_1)\cap
\mathfrak{F}_P(A_2)\bigr)\ .
\end{equation}

Let
$
\bmap:\Ant(P)\rightarrow\Ant(P)
$
be the {\em blocker map on $\Ant(P)$}; by definition,
\begin{equation}
\bmap:A\mapsto\bmin\mathbf{I}(P,A)\ .
\end{equation}
In particular, for every $a\in P$, $a>_P\hat{0}_P$, we have
$\bmap(\{a\})=\mathfrak{I}_P(a)\cap\Pa$. We also have
\begin{equation}
\bmap(\emptyset\subset P)=\{\hat{0}_P\}\ ,\ \ \
\bmap\bigl(\{\hat{0}_P\}\bigr)=\emptyset\subset P\ .
\end{equation}

For a one-element antichain $\{a\}$, we write $\bmap(a)$ instead
of $\bmap(\{a\})$.

If $A$ is a nontrivial antichain of $P$ then
Lemma~\ref{Intersectersandcomplementers4} implicitly states the
following equalities in $\mathfrak{A}(P)$:
\begin{equation}
\label{Blockermapandcomplementarymap3}
\begin{split}
\bmap(A)=\bigwedge_{a\in A}\quad\bigvee_{e\in\bmap(a)} \{e\}=
\bigvee_{E\in\mathcal{B}(\{\bmap(a):a\in A\})}
\quad\bigwedge_{e\in E}\{e\}\ .
\end{split}
\end{equation}

Let $\Antb(P)$ denote the image of $\Ant(P)$ under the blocker
map. The set $\Antb(P)$ is equipped, by definition, with the
partial order induced by the partial order on $\Ant(P)$. For a
blocker $B\in\Antb(P)$, the subposet
$\bmap^{-1}(B)=\{A\in\Ant(P):\bmap(A)=B\}$ is the preimage of $B$
under the blocker map.

The following lemma is a reformulation of
~(\ref{Intersectersandcomplementers3}):
\begin{lem}
\label{Blockermapandcomplementarymap1} If $A_1,A_2\in\Ant(P)$ and
$A_1\leq_{\Ant(P)} A_2$ then $\bmap(A_1)\geq_{\Antb(P)}\bmap(A_2)$.
\end{lem}

Definition~\ref{Intersectersandcomplementers1} implies the following {\em
reciprocity property for intersecters}: for every antichain $A$ of $P$, we have
\begin{equation}
A\subseteq\mathbf{I}\bigl(P,\bmap(A)\bigr)\ .
\end{equation}

In the theory of blocking sets the following fact is basic:
\begin{prop}{\rm(see~\cite{EF,L})}.
For any clutter $\mathcal{G}$,
$\mathcal{B}\bigl(\mathcal{B}(\mathcal{G})\bigr)=\mathcal{G} $.
\end{prop}

This statement may be generalized in the following way:
\begin{thm}
\label{Blockermapandcomplementarymap2} The restriction map
$\bmap\!\!\mid_{\Antb(P)}$ is an involution, that is, for each
blocker $B\in\Antb(P)$, $\bmap\bigl(\bmap(B)\bigr)=B$.
\end{thm}

\begin{proof}
There is nothing to prove for the {\em trivial blockers\/}
$B=\hat{0}_{\Antb(P)}=\emptyset\subset P$ and
$B=\hat{1}_{\Antb(P)}=\{\hat{0}_P\}$. So suppose that $B$ is
nontrivial. Choose an arbitrary antichain $A'\in\bmap^{-1}(B)$.
With regard to reciprocity property for intersecters, every
element of $A'$ is an intersecter for the antichain $B=\bmap(A')$.
In other words, for each element $a'\in A'$ we have the inclusion
$a'\in\mathbf{I}(P,B)= \bigcap_{b\in
B}\mathfrak{F}_P\bigl(\bmap(b)\bigr)$. Taking this inclusion into
account, we assign to the antichain $A'$ the antichain
\begin{equation}
A=\bmin\bigcap_{b\in B}\mathfrak{F}_P\bigl(\bmap(b)\bigr)\in\bmap^{-1}(B)\ ,
\end{equation}
which is the blocker of $B$, by (\ref{eq:2006-1}). Then
$\bmap(A)=B,\ \bmap(B)=A$, and the theorem follows.
\end{proof}

By Lemma~\ref{Intersectersandcomplementers4}, a nontrivial
antichain $A$ of $P$, considered as an element of $\Ant(P)$, is a
fixed point of the blocker map on $\Ant(P)$ if and only if $A=
\bigwedge_{a\in A}\bigvee_{e\in\bmap(a)} \{e\}$ or, equivalently,
$A= \bigvee_{E\in\mathcal{B}(\{\bmap(a):a\in A\})}\bigwedge_{e\in
E}\{e\}$. We study the structure of a preimage of the blocker map.

\begin{thm}
For each blocker $B\in\Antb(P)$, its preimage $\bmap^{-1}(B)$ is a
join-subsemilattice of the lattice $\Ant(P)$.
\end{thm}

\begin{proof}
There is nothing to prove for a trivial blocker $B$, so suppose that $B$ is
nontrivial. Choose two antichains $A_1,A_2\in\bmap^{-1}(B)$. According
to~(\ref{Blockermapandcomplementarymap3}), we have the following equalities in
the lattice $\Ant(P)$:
\begin{equation}
B=\bmap(A_1)=\bigwedge_{a_1\in A_1}\quad\bigvee_{e\in\bmap(a_1)}
\{e\}=\bmap(A_2)=\bigwedge_{a_2\in A_2}\quad\bigvee_{e\in\bmap(a_2)} \{e\}\ .
\end{equation}
Therefore
\begin{equation}
B=\bigwedge_{a\in A_1\vee_{\Ant(P)}A_2}\quad\bigvee_{e\in\bmap(a)}
\{e\}=\bmap(A_1\vee_{\Ant(P)}A_2)\ .
\end{equation}
Hence $A_1\vee_{\Ant(P)}A_2\in\bmap^{-1}(B)$.
\end{proof}

The greatest element of $\bmap^{-1}(B)$ is $\bmap(B)$.

Let $\cmap:\Ant(P)\rightarrow\Ant(P)$ be the {\em complementary map on
$\Ant(P)$}; by definition,
\begin{equation}
\cmap:A\mapsto\bmax\mathbf{C}(P,A)\ .
\end{equation}
In particular, we have $\cmap(\emptyset\subset P)=\emptyset\subset P$ and
$\cmap(\{\hat{0}_P\})=\{\hat{1}_P\}$.

Let $\Antc(P)$ denote the image of $\Ant(P)$ under the
complementary map. The set $\Antc(P)$ is equipped, by definition,
with the partial order induced by the partial order on the lattice
of order ideals of $P$: for $C_1,C_2\in\Antc(P)$ we set
$C_1\leq_{\Antc(P)}C_2$ if and only if
$\mathfrak{I}_P(C_1)\subseteq\mathfrak{I}_P(C_2)$.

\section{Lattice of blockers}
\label{Latticeofblockers}

In this section, we study the structure of the poset of blockers
in $P$.

\begin{lem}
The poset $\Antb(P)$ of blockers in $P$ is a meet-subsemilattice of the lattice
$\Ant(P)$.
\end{lem}

\begin{proof}
We have to prove that for all $B_1,B_2\in\Antb(P)$, it holds
$B_1\wedge_{\Ant(P)}B_2\in\Antb(P)$. There is nothing to prove when one of the
blockers $B_1,B_2$ is trivial. Suppose that both $B_1$ and $B_2$ are
nontrivial. With the help of Theorem~\ref{Blockermapandcomplementarymap2}, we
write
\begin{equation}
B_1\wedge_{\Ant(P)}B_2= \bmap\bigl(\bmap(B_1)\bigr) \wedge_{\Ant(P)}
\bmap\bigl(\bmap(B_2)\bigr)\ .
\end{equation}
According to~(\ref{Blockermapandcomplementarymap3}), we have the following
equalities in $\Ant(P)$:
\begin{equation}
\begin{split}
B_1\wedge_{\Ant(P)}B_2&=
\left(\bigwedge_{a_1\in\bmap(B_1)}\quad\bigvee_{e\in\bmap(a_1)} \{e\}\right)
\wedge_{\Ant(P)}
\left(\bigwedge_{a_2\in\bmap(B_2)}\quad\bigvee_{e\in\bmap(a_2)}
\{e\}\right)\\&=
\bigwedge_{a\in\bmap(B_1)\vee_{\Ant(P)}\bmap(B_2)}\quad\bigvee_{e\in\bmap(a)}
\{e\}= \bmap\bigl( \bmap(B_1) \vee_{\Ant(P)} \bmap(B_2) \bigr)\in\Antb(P)\ .
\end{split}
\end{equation}
\end{proof}

\begin{lem}
The meet-semilattice $\Antb(P)$ is self-dual.
\end{lem}

\begin{proof}
Let $B_1,B_2\in\Antb(P)$. If $B_1\leq_{\Antb(P)} B_2$ then
$B_1\leq_{\Ant(P)}B_2$, and we see that $\bmap(B_1)\geq_{\Antb(P)}\bmap(B_2)$,
by Lemma~\ref{Blockermapandcomplementarymap1}.

Conversely, the relation $\bmap(B_1)\geq_{\Antb(P)}\bmap(B_2)$ implies the
relation $B_1=\bmap\bigl(\bmap(B_1)\bigr)\leq_{\Antb(P)}
B_2=\bmap\bigl(\bmap(B_2)\bigr)$, in view of
Theorem~\ref{Blockermapandcomplementarymap2} and
Lemma~\ref{Blockermapandcomplementarymap1}.

Because the restriction map $\bmap\!\!\mid_{\Antb(P)}$ is bijective, we see
that it is an antiautomorphism of $\Antb(P)$.
\end{proof}

We now summarize the information of this section:

\begin{thm}
\label{Latticeofblockers3} The poset $\Antb(P)$ is a lattice with the least
element $\hat{0}_{\Antb(P)}=\emptyset\subset P$ and greatest element
$\hat{1}_{\Antb(P)}=\{\hat{0}_P\}$. The unique atom of $\Antb(P)$ is
$\bmap(\Pa)$, and the unique coatom of $\Antb(P)$ is $\Pa$. Moreover,
\begin{enumerate}
\renewcommand{\theenumi}{\rm{\roman{enumi}}}

\item
the poset $\Antb(P)$ is a meet-subsemilattice of the lattice
$\Ant(P)$,

\item
\label{Latticeofblockers4} the lattice $\Antb(P)$ is self-dual,

\item
in the lattice $\Antb(P)$ the operations of meet and join are
determined as follows: for $B_1,B_2\in\Antb(P)$
\begin{gather}
\label{Latticeofblockers1} B_1\wedge_{\Antb(P)}B_2=B_1\wedge_{\Ant(P)}B_2\\
\intertext{and} \label{Latticeofblockers2} B_1\vee_{\Antb(P)}B_2=
\bmap\bigl(\bmap(B_1)\wedge_{\Ant(P)}\bmap(B_2)\bigr)\ .
\end{gather}
\end{enumerate}
\end{thm}

\begin{proof}
The only missing step is to prove~{\rm(\ref{Latticeofblockers2})},
but the equality $B_1\vee_{\Antb(P)}B_2=
\bmap\bigl(\bmap(B_1)\wedge_{\Antb(P)}\bmap(B_2)\bigr)$
immediately follows from the self-duality of the lattice
$\Antb(P)$, in view of the existence of its anti-automorphism
$\bmap\!\!\mid_{\Antb(P)}$. With the help of
equality~{\rm(\ref{Latticeofblockers1})}, we
obtain~{\rm(\ref{Latticeofblockers2})}.
\end{proof}

We call the lattice $\Antb(P)$ the {\em lattice of blockers in the
poset $P$}. It follows immediately from the definition of the
complementary map that its restriction $\cmap\!\!\mid_{\Antb(P)}:
\Antb(P)\rightarrow\Antc(P)$, $B\mapsto\cmap(B)$, is an
isomorphism of $\Antb(P)$ into the lattice $\Antc(P)$.


\end{document}